\documentclass[12pt,a4paper]{article}
\usepackage{amsmath}
\usepackage{amssymb}
\usepackage{t1enc}
\usepackage[latin1]{inputenc}
\usepackage[english]{babel}
\usepackage{amsfonts}
\usepackage{latexsym}
\usepackage{bm}
\newtheorem{theorem}{Theorem}[section]

\newtheorem{lemma}[theorem]{Lemma}

\newtheorem{conj}[theorem]{Conjecture}
\newtheorem{remark}[theorem]{Remark}

\title{Cross-intersecting subfamilies of levels of hereditary families}
\author{Peter Borg\footnote{Address: Department of Mathematics, University of Malta, Malta. \newline Email: peter.borg@um.edu.mt}
}
\date{}

\begin{document}
\maketitle

\begin{abstract}
A set $A$ \emph{$t$-intersects} a set $B$ if $A$ and $B$ have at least $t$ common elements. Families $\mathcal{A}_1, \mathcal{A}_2, \dots, \mathcal{A}_k$ of sets are \emph{cross-$t$-intersecting} if, for every $i$ and $j$ in $\{1, 2, \dots, k\}$ with $i \neq j$, each set in $\mathcal{A}_i$ $t$-intersects each set in $\mathcal{A}_j$. An active problem in extremal set theory is to determine, for a given finite family $\mathcal{F}$, the structure of $k$ cross-$t$-intersecting subfamilies whose sum or product of sizes is maximum. For a family $\mathcal{H}$, the \emph{$r$-th level} $\mathcal{H}^{(r)}$ of $\mathcal{H}$ is the family of all sets in $\mathcal{H}$ of size $r$, and, for $s \leq r$, $\mathcal{H}^{(s)}$ is called a \emph{$(\leq r)$-level of $\mathcal{H}$}. We solve the problem for any union $\mathcal{F}$ of $(\leq r)$-levels of any union $\mathcal{H}$ of power sets of sets of size at least a certain integer $n_0$, where $n_0$ is independent of $\mathcal{H}$ and $k$ but depends on $r$ and $t$ (dependence on $r$ is inevitable, but dependence on $t$ can be avoided). Our primary result asserts that there are only two possible optimal configurations for the sum. %and characterizes them. A special case of our result for the sum
A special case was conjectured by Kamat in 2011.
We also prove generalizations, whereby $\mathcal{A}_1, \mathcal{A}_2, \dots, \mathcal{A}_k$ are not necessarily contained in the same union of levels. Various Erd\H os--Ko--Rado-type results follow. %, particularly the result in [J. London Math. Soc. 79 (2009), 167--185].
The sum problem for a level of a power set was solved for $t=1$ by Hilton in 1977, %(the first result of this kind), 
and for any $t$ by Wang and Zhang in 2011.
\end{abstract}

\section{Introduction} \label{Intro}

Before introducing the problems treated in this paper, we provide the main definitions and notation.

Unless otherwise stated, we shall use small letters such as $x$ to denote non-negative integers or elements of a set, capital letters such as $X$ to denote sets, and calligraphic letters such as $\mathcal{F}$ to denote \emph{families} (that is, sets whose members are sets themselves). Arbitrary sets and families are taken to be finite and may be the \emph{empty set} $\emptyset$. An \emph{$r$-set} is a set of size $r$, that is, a set having exactly $r$ elements (also called members).

The set $\{1, 2, \dots\}$ of positive integers is denoted by
$\mathbb{N}$. For $m, n \in \mathbb{N}$ with $m \leq n$, the set $\{i \in \mathbb{N} \colon m \leq i \leq n\}$ is denoted by $[m,n]$, and $[1,n]$ is abbreviated to $[n]$. For a set $X$, the \emph{power set of $X$} (that is, $\{A \colon A \subseteq X\}$) is denoted by $2^X$, and the families $\{Y \subseteq X \colon |Y| = r\}$ and $\{Y \subseteq X\colon |Y| \leq r\}$ are denoted by $X \choose r$ and $X \choose \leq r$, respectively.

For any family $\mathcal{F}$, we have the following definitions and notation. A \emph{base of $\mathcal{F}$} is a set $B$ in $\mathcal{F}$ such that, for each $A$ in $\mathcal{F}$, $B$ is not a proper subset of $A$. %is a set in $\mathcal{F}$ that is a subset of no other set in $\mathcal{F}$. 
The size of a smallest base of $\mathcal{F}$ is denoted by $\mu(\mathcal{F})$. %The union of all sets in $\mathcal{F}$ is denoted by $U(\mathcal{F})$.
The families $\{F \in \mathcal{F} \colon |F| = r\}$ and $\{F \in \mathcal{F} \colon |F| \leq r\}$ are denoted by $\mathcal{F}^{(r)}$ and $\mathcal{F}^{(\leq r)}$, respectively. The family $\mathcal{F}^{(r)}$ is called the \emph{$r$-th level of $\mathcal{F}$}, and, for $s \leq r$, $\mathcal{F}^{(s)}$ is called a \emph{$(\leq r)$-level of $\mathcal{F}$}. For any set $T$, we denote $\{F \in \mathcal{F} \colon T \subseteq F\}$ by $\mathcal{F}\langle T \rangle$. We call $\mathcal{F}\langle T \rangle$ a \emph{$t$-star of $\mathcal{F}$} if $|T|=t$. %and $\mathcal{F}\langle T \rangle \neq \emptyset$.

Given an integer $t \geq 1$, we say that a set $A$
\emph{$t$-intersects} a set $B$ if $A$ and $B$ have at least $t$
common elements. A family $\mathcal{A}$ is said to be
\emph{$t$-intersecting} if, for every $A, B \in \mathcal{A}$, $A$ $t$-intersects $B$. %each set in $\mathcal{A}$ $t$-intersects all the other sets in $\mathcal{A}$ (that is, $|A \cap B| \geq t$ for every $A, B \in \mathcal{A}$ with $A \neq B$). 
A $1$-intersecting family is also simply called an
\emph{intersecting family}. A $t$-intersecting family $\mathcal{A}$ is said to be \emph{trivial} if its sets have at least $t$ common elements (that is, $\left| \bigcap_{A \in \mathcal{A}} A \right| \geq t$). Note that non-empty $t$-stars are trivial $t$-intersecting families. We say that a family $\mathcal{F}$ has the \emph{$t$-star property} if at least one of the largest $t$-intersecting subfamilies of $\mathcal{F}$ is a $t$-star of $\mathcal{F}$. We also say that $\mathcal{F}$ has the \emph{strict $t$-star property} if all the largest $t$-intersecting subfamilies of $\mathcal{F}$ are $t$-stars of $\mathcal{F}$.

The study of intersecting families started in \cite{EKR}, which features the classical result, known as the Erd\H os--Ko--Rado (EKR) Theorem, that says that, for $1 \leq t \leq r$, there exists an integer $n_0(r,t)$ such that, for every $n \geq n_0(r,t)$, the size of a largest $t$-intersecting subfamily of ${[n] \choose r}$ is the size ${n-t \choose r-t}$ of every non-empty $t$-star of ${[n] \choose r}$, meaning that ${[n] \choose r}$ has the $t$-star property. It was also shown in \cite{EKR} that the smallest possible value of $n_0(r,1)$ is $2r$; %and two of the various proofs of this fact (see \cite{Kat,K,D}) are particularly short and beautiful: Katona's \cite{K}, introducing the elegant cycle method, and Daykin's \cite{D}, using the powerful Kruskal--Katona Theorem \cite{Ka,Kr}. 
among the various proofs of this fact (see \cite{EKR,Kat,HM,K,D,FF2,HK2}) 
there is a short one by Katona \cite{K}, introducing the elegant cycle method, and another one by Daykin \cite{D}, using the Kruskal-Katona Theorem \cite{Ka,Kr}. Note that, for $n/2 < r < n$, ${[n] \choose r}$ itself is intersecting and hence does not have the $1$-star property. A sequence of results \cite{EKR,F_t1,W,FF,AK1} culminated in the complete solution of the problem for $t$-intersecting subfamilies of ${[n] \choose r}$. The solution confirmed a conjecture of Frankl \cite{F_t1}. Frankl \cite{F_t1} and Wilson \cite{W} proved the following.
\begin{theorem}[\cite{F_t1,W}]\label{AK} Let $1 \leq t < r < n$. Then: \\
(i) ${[n] \choose r}$ has the $t$-star property if and only if $n \geq (r-t+1)(t+1)$.\\
(ii) ${[n] \choose r}$ has the strict $t$-star property if and only if $n > (r-t+1)(t+1)$.
\end{theorem}
Trivially, for $t = r$ or $r = n$, a $t$-intersecting subfamily of ${[n] \choose r}$ can have only one member, so ${[n] \choose r}$ has the strict $t$-star property. The $t$-intersection problem for $2^{[n]}$ was solved by Katona \cite{Kat}. These are among the most prominent results in extremal set theory. The EKR Theorem inspired a wealth of results that establish how large a system of sets can be under certain intersection conditions; see \cite{DF,F,F2,Borg7,HST,HT,FTsurvey}.

A family $\mathcal{F}$ is said to be \emph{hereditary} if, for each $F \in \mathcal{F}$, all the subsets of $F$ are members of $\mathcal{F}$. The power set is the simplest example. In fact, by definition, a family is hereditary if and only if it is a union of power sets. Note that, if $X_1, \dots, X_k$ are the bases of a hereditary family $\mathcal{H}$, then $\mathcal{H} = 2^{X_1} \cup \dots \cup 2^{X_k}$.

Hereditary families are important combinatorial objects that have attracted much attention. In the literature, a hereditary family is also called an \emph{ideal}, a \emph{downset}, and an \emph{abstract simplicial complex}. The various interesting examples include the family of \emph{independent sets} of a \emph{graph} or of a \emph{matroid}. One of the central problems in extremal set theory is Chv\'atal's conjecture \cite{Chv}, which claims that every hereditary family $\mathcal{H}$ has the $1$-star property. The best result so far on this conjecture is due to Snevily \cite{Sn} (see also \cite{Borg4b}). %provides a generalization obtained by means of a self-contained alternative argument). 
The conjecture cannot be generalized for $t$-intersecting subfamilies. Indeed, if $2 \leq t < n$ and $\mathcal{H} = 2^{[n]}$, then $\mathcal{H}$ does not have the $t$-star property; the complete characterization of the largest $t$-intersecting subfamilies of $2^{[n]}$ is given in \cite{Kat}. For levels of hereditary families, we have the following generalization of the Holroyd--Talbot Conjecture \cite[Conjecture~7]{HT}.
\begin{conj}[\cite{Borg}] \label{AK gen} If $1 \leq t \leq r$, $\emptyset \neq  S \subseteq [t,r]$, and $\mathcal{H}$ is a hereditary family with $\mu(\mathcal{H}) \geq (t+1)(r-t+1)$, then \\
(i) $\bigcup_{s \in S}\mathcal{H}^{(s)}$ has the $t$-star property, \\
(ii) $\bigcup_{s \in S}\mathcal{H}^{(s)}$ has the strict $t$-star property if $\mu(\mathcal{H}) > (t+1)(r-t+1)$ or $S \neq \{r\}$.
\end{conj}
Note that, if $\mathcal{H} = 2^{[n]}$, then $\mathcal{H}^{(s)} = {[n] \choose s}$ and $\mu(\mathcal{H}) = n$. Clearly, it follows by
Theorem~\ref{AK} that the conjecture is true for
$\mathcal{H}=2^{[n]}$ and that the condition $\mu(\mathcal{H}) \geq (t+1)(r-t+1)$ cannot be improved (in fact, one can check that, for $t < r < n < (t+1)(r-t+1)$, we have $|\{A \in {[n] \choose r} \colon |A \cap [t+2]| \geq t+1\}| > {n-t \choose r-t}$, and so the $t$-stars of ${[n] \choose r}$ are not among the largest $t$-intersecting subfamilies of ${[n] \choose r}$). The conjecture is true for $\mu(\mathcal{H})$ sufficiently large depending only on $r$ and $t$.
\begin{theorem}[\cite{Borg}]\label{t int her} Conjecture~\ref{AK gen} is true if $\mu(\mathcal{H}) \geq (r-t){3r-2t-1 \choose t+1} + r$.
\end{theorem}
\begin{remark}\label{rem1} \emph{By weakening the condition on $\mu(\mathcal{H})$ in Theorem~\ref{t int her}, we can eliminate the dependence on $t$. Indeed, we crudely have $(r-t){3r-2t-1 \choose t+1} + r < r{3r \choose t+1} < r{3r \choose \lfloor 3r/2 \rfloor}$. Thus, Conjecture~\ref{AK gen} is true if $\mu(\mathcal{H}) \geq r{3r \choose \lfloor 3r/2 \rfloor}$. Dependence on $r$ is inevitable, given that, as pointed out above, the condition $\mu(\mathcal{H}) \geq (t+1)(r-t+1)$ in Conjecture~\ref{AK gen} cannot be improved.}
\end{remark}

A problem that generalizes the intersection problem described above and that is also attracting much attention is the cross-intersection problem.

Families $\mathcal{A}_1, \dots, \mathcal{A}_k$ are said to be \emph{cross-$t$-intersecting} if, for every $i$ and $j$ in $[k]$ with $i \neq j$, each set in $\mathcal{A}_i$ $t$-intersects each set in $\mathcal{A}_j$. Cross-$1$-intersecting families are also simply called \emph{cross-intersecting families}. Note that, if $\mathcal{A}$ is a $t$-intersecting family and $\mathcal{A}_1 = \dots = \mathcal{A}_k = \mathcal{A}$, then $\mathcal{A}_1, \dots, \mathcal{A}_k$ are cross-$t$-intersecting.

For $t$-intersecting subfamilies of a given family $\mathcal{F}$,
the natural question to ask is how large they can be. For
cross-$t$-intersecting families, two natural parameters arise: the
sum and the product of sizes of the cross-$t$-intersecting
families (note that the product of sizes of $k$ families
$\mathcal{A}_1, \dots, \mathcal{A}_k$ is the number of $k$-tuples
$(A_1, \dots, A_k)$ such that $A_i \in \mathcal{A}_i$ for each $i
\in [k]$). It is therefore natural to consider the problem of
maximizing the sum or the product of sizes of $k$
cross-$t$-intersecting subfamilies (not necessarily distinct or
non-empty) of a given family $\mathcal{F}$. The paper \cite{Borg8} analyses this problem in general and shows in particular that, for $k$ sufficiently large, both the sum and the product are maxima if $\mathcal{A}_1 = \dots = \mathcal{A}_k = \mathcal{L}$ for some largest $t$-intersecting subfamily $\mathcal{L}$ of $\mathcal{F}$. Therefore, this problem incorporates the $t$-intersection problem.  Solutions have been obtained for various families; many results are outlined in \cite{Borg8}. In this paper we are primarily concerned with the case where $\mathcal{F}$ is a level or a union of levels of a hereditary family, as in Theorem~\ref{t int her}, but we also consider even more general settings. Before coming to the contributions in this paper, we shall outline more facts and provide further motivation.

The cross-intersection problem described above has its origin in \cite{H}, in which Hilton solved the sum problem for ${[n] \choose r}$ and $t = 1$. Wang and Zhang \cite{WZ} solved the sum problem for ${[n] \choose r}$ and any $t$ by reducing it to the complete $t$-intersection theorem of Ahlswede and Khachatrian \cite{AK1}, using a striking combination of the method in \cite{Borg3,Borg2,BL2,Borg5} and the \emph{no-homomorphism lemma} \cite{AC,CK}. Their result holds for a much more general class of important families and was used in the solution of the sum problem for $2^{[n]}$ \cite[Theorems~3.10, 4.1]{Borg8}.

The maximum product problem for $2^{[n]}$ was settled in \cite{MT2} for the case where $k = 2$ or $n+t$ is even (see \cite[Section~5.2]{Borg8}, which features a conjecture for the case where $k > 2$ and $n+t$ is odd). As pointed out above, ${[n] \choose r}$ is the $r$-th level of $2^{[n]}$. Pyber \cite{Pyber} proved that, for any $r$, $s$ and $n$ such that either $1 \leq s = r \leq n/2$ or $1 \leq s < r$ and $n \geq 2r + s -2$, if $\mathcal{A} \subseteq {[n] \choose r}$ and $\mathcal{B} \subseteq {[n] \choose s}$ such that $\mathcal{A}$ and $\mathcal{B}$ are cross-intersecting, then $|\mathcal{A}||\mathcal{B}| \leq {n-1 \choose r-1}{n-1 \choose s-1}$. Subsequently, Matsumoto and Tokushige \cite{MT} proved this for $1 \leq s \leq r \leq n/2$. As \cite[Conjecture~3]{Hirschorn} states, if $\mathcal{A}$ and $\mathcal{B}$ are cross-$t$-intersecting and $n$ is sufficiently large, then $|\mathcal{A}||\mathcal{B}| \leq {n-t \choose r-t}{n-t \choose s-t}$; see Theorem~\ref{prodcor2}, which covers the general case with $k \geq 2$ cross-$t$-intersecting families. 

This brings us to the results in this paper, which we present formally in the next section. We mainly determine the optimal structures for both the sum problem and the product problem for unions of levels of hereditary families with sufficiently large bases. We give an affirmative answer to a generalization of a conjecture of Kamat \cite{Kamat} regarding the sum problem for the families just mentioned. One important aspect of each of our main results is that the condition for how large the bases of a hereditary family $\mathcal{H}$ should be is independent of $\mathcal{H}$ and the number $k$ of cross-$t$-intersecting families, but depends only on $t$ and the maximum possible size $r$ of the sets in the levels (as in Conjecture~\ref{AK gen}), and can also be made independent of $t$ as in Remark~\ref{rem1} (dependence on $r$ is inevitable). 

Hereditary families exhibit undesirable phenomena. The motivation behind establishing a cross-$t$-intersection result for a \emph{union} of levels (similarly to Theorem~\ref{t int her}) is that, for a hereditary family, this general form cannot be immediately deduced from the result for just one level; see \cite[Example~1]{Borg}. The complete absence of symmetry makes intersection problems like the ones described above difficult to deal with. Many of the well-known techniques in extremal set theory, such as the EKR shifting technique (see \cite{EKR,F}) and Katona's cycle method \cite{K}, fail to work for hereditary families. The ingredients that enable us to overcome such difficulties (for the problems we address in this paper) are given in Sections~\ref{keysection} and \ref{propertysection}.

\section{Results}

This section is divided into two subsections. In the first one we present the above-mentioned general results for the maximum sum, and in the second we present those for the maximum product. We also point out various consequences.

Throughout this section and the rest of the paper, we take
\begin{gather} m(r,t) := \max \left\{ 2r, \; \frac{(r-t)(r-t+5)}{2} +
(t-1) \right\}, \nonumber \\
n_{\rm S}(r,t) := (r-t+1){m(r,t) \choose t+1} + r, \nonumber \\
n_{\rm P}(r,t) := (r-t){r \choose t}{m(r,t) \choose t+1} + r. \nonumber
\end{gather}

\subsection{The maximum sum} \label{maxsumsubsec}

%We now present various cross-$t$-intersection results that concern families of particular importance and that follow from Theorem~\ref{sumresult}. 
The following is our result for the sum cross-intersection problem described in Section~\ref{Intro}.

\begin{theorem} \label{sumcor2} If $1 \leq t \leq r$, $\mathcal{H}$ is a hereditary family with $\mu(\mathcal{H}) \geq n_{\rm S}(r,t)$, $\emptyset \neq S \subseteq [t,r]$, $\mathcal{F} := \bigcup_{s \in S}\mathcal{H}^{(s)}$, $\mathcal{A}_1, \dots, \mathcal{A}_k$ are cross-$t$-intersecting subfamilies of $\mathcal{F}$, and $\mathcal{L}$ is a largest $t$-star of $\mathcal{F}$, then
\[\sum_{i = 1}^k |\mathcal{A}_i| \leq \max\left\{ k
|\mathcal{L}|, |\mathcal{F}| \right\}.\]
Moreover, equality holds if and only if one of the
following holds: \\
(i) $|\mathcal{F}| < k|\mathcal{L}|$ and $\mathcal{A}_1 = \dots = \mathcal{A}_k = \mathcal{L}'$ for some largest $t$-star
$\mathcal{L}'$ of $\mathcal{F}$.\\
(ii) $|\mathcal{F}| > k |\mathcal{L}|$ and, for some $j \in [k]$,
$\mathcal{A}_j = \mathcal{F}$ and $\mathcal{A}_i = \emptyset$ for
all $i \in [k] \backslash \{j\}$.\\
(iii) $|\mathcal{F}| = k |\mathcal{L}|$, and $\mathcal{A}_1,
\dots, \mathcal{A}_k$ are as in (i) or (ii).
\end{theorem}
We will show that this follows immediately from Theorem~\ref{sumresult}. Note that, as in Remark~\ref{rem1}, we can obtain a condition for $\mu(\mathcal{H})$ that is dependent only on $r$. 

\begin{remark} \emph{Theorem~\ref{sumcor2} implies Theorem~\ref{t int her} for $\mu(\mathcal{H}) \geq n_{\rm S}(r,t)$. Indeed, if $k > |\mathcal{F}|/|\mathcal{L}|$, $\mathcal{A}$ is an intersecting subfamily of $\mathcal{F}$, and $\mathcal{A}_1 = \dots = \mathcal{A}_k = \mathcal{A}$, then (i) holds, $\mathcal{A}_1, \dots, \mathcal{A}_k$ are cross-$t$-intersecting subfamilies of $\mathcal{F}$, and hence $|\mathcal{A}| = \frac{1}{k} \sum_{i=1}^k |\mathcal{A}_i| \leq |\mathcal{L}|$.}
\end{remark}

If one of $k$ cross-$t$-intersecting families has a member of size less than $t$, then the other $k-1$ families are empty. Consequently, as indicated in the proof of Theorem~\ref{sumcor2.5}, Theorem~\ref{sumcor2} and the subsequent results for $\sum_{i = 1}^k |\mathcal{A}_i|$ and $\prod_{i = 1}^k |\mathcal{A}_i|$ actually hold for any union of $(\leq r)$-levels of $\mathcal{H}$ (that is, we can take $\emptyset \neq S \subseteq \{0\} \cup [r]$ rather than $\emptyset \neq S \subseteq [t,r]$). 

Of particular interest is the case where $\mathcal{F}$ consists of all the sets in $\mathcal{H}$ of size at most $r$, in which case $\mathcal{F}$ itself is hereditary.

\begin{theorem} \label{sumcor2.5} If $1 \leq t \leq r$, $\mathcal{H}$ is a hereditary family with $\mu(\mathcal{H}) \geq n_{\rm S}(r,t)$, $\mathcal{A}_1, \dots, \mathcal{A}_k$ are cross-$t$-intersecting subfamilies of $\mathcal{H}^{(\leq r)}$, and $\mathcal{L}$ is a largest $t$-star of $\mathcal{H}^{(\leq r)}$, then
\[\sum_{i = 1}^k |\mathcal{A}_i| \leq \max\left\{ k
|\mathcal{L}|, |\mathcal{H}^{(\leq r)}| \right\}.\]
%
%Moreover, equality holds if and only of one of (i)--(iii) of Theorem~\ref{sumcor2} holds with $\mathcal{H}^{(\leq r)}$ instead of $\mathcal{F}$.
\end{theorem}
\textbf{Proof.} Let $\mathcal{A}_1, \dots, \mathcal{A}_k$ be such that $\sum_{i = 1}^k |\mathcal{A}_i|$ is maximum. Let $I = \{i \in [k] \colon \mathcal{A}_i \neq \emptyset\}$. Then $|I| \geq 1$. Let $h \in I$. If $I = \{h\}$, then $\mathcal{A}_h = \mathcal{H}^{(\leq r)}$ and $\mathcal{A}_i = \emptyset$ for all $i \in [k] \backslash \{h\}$. Suppose $I \neq \{h\}$. Let $j \in I \backslash \{h\}$. Let $A_1 \in \mathcal{A}_h$ and $A_2 \in \mathcal{A}_j$. For each $i \in [k] \backslash \{h\}$, $|A| \geq |A \cap A_1| \geq t$ for any $A \in \mathcal{A}_i$. For each $A \in \mathcal{A}_h$, $|A| \geq |A \cap A_2| \geq t$.  Thus, $\mathcal{A}_1, \dots, \mathcal{A}_k \subseteq \bigcup_{s \in [t,r]}\mathcal{H}^{(s)}$. Since $|I| \geq 2$, it follows by the choice of $\mathcal{A}_1, \dots, \mathcal{A}_k$ and by Theorem~\ref{sumcor2} that $\mathcal{A}_1 = \dots = \mathcal{A}_k = \mathcal{L}'$ for some largest $t$-star $\mathcal{L}'$ of $\bigcup_{s \in [t,r]}\mathcal{H}^{(s)}$.~\hfill{$\Box$} \\

%If there exists some $h \in [k]$ such that $\mathcal{A}_h$ has a member of size at most $t-1$, then, since $\mathcal{A}_1, \dots, \mathcal{A}_k$ are cross-$t$-intersecting, $\mathcal{A}_j = \emptyset$ for all $j \in [k] \backslash \{h\}$, so $\sum_{i = 1}^k |\mathcal{A}_i| = |\mathcal{A}_h| \leq |\mathcal{H}^{(\leq r)}|$. If none of $\mathcal{A}_1, \dots, \mathcal{A}_k$ has a member of size at most $t-1$, then $\mathcal{A}_1, \dots, \mathcal{A}_k \subseteq \bigcup_{s \in [t,r]}\mathcal{H}^{(s)}$, so the result follows by Theorem~\ref{sumcor2}.~\hfill{$\Box$} \\

We have actually shown that the optimal structures for Theorem~\ref{sumcor2.5} are given by Theorem~\ref{sumcor2} (i)--(iii) with $\mathcal{H}^{(\leq r)}$ instead of $\mathcal{F}$.

The next theorem is our most general result for the maximum sum. It not only allows each family $\mathcal{A}_i$ %(of the $k$ cross-$t$-intersecting families) 
to be contained in an arbitrary union $\mathcal{F}_i$ of $(\leq r)$-levels, but also allows the families $\mathcal{F}_1, \dots, \mathcal{F}_k$ to be different. It is proved in Section~\ref{sumsection}.

\begin{theorem} \label{sumresult} Let $1 \leq t \leq r$, and let $\mathcal{H}$ be a hereditary family with $\mu(\mathcal{H}) \geq n_{\rm S}(r,t)$. Let $k \geq 2$, and, for each $i \in [k]$, let $\emptyset \neq S_i \subseteq [t,r]$, $\mathcal{F}_i := \bigcup_{s \in
S_i}\mathcal{H}^{(s)}$ and $\mathcal{T}_i := \{T \in \mathcal{H}^{(t)} \colon \mathcal{F}_i \langle T \rangle$ is a largest
$t$-star of $\mathcal{F}_i\}$. Let $\mathcal{A}_1, \dots, \mathcal{A}_k$ be cross-$t$-intersecting families such that $\mathcal{A}_i
\subseteq \mathcal{F}_i$ for all $i \in [k]$. Then $\mathcal{T}_i \neq \emptyset$ for all $i \in [k]$, and, for any $T_1 \in \mathcal{T}_1, \dots, T_k \in \mathcal{T}_k$,
\[\sum_{i = 1}^k |\mathcal{A}_i| \leq \max\left\{\sum_{i = 1}^k
|\mathcal{F}_i \langle T_i \rangle|, |\mathcal{F}_1|, \dots, |\mathcal{F}_k|\right\}.\]
Moreover, equality holds if and only if one of the following
holds: \\
(i) $\sum_{i = 1}^k |\mathcal{F}_i \langle T_i \rangle| >
\max\{|\mathcal{F}_1|, \dots, |\mathcal{F}_k|\}$,
$\bigcap_{i = 1}^k \mathcal{T}_i \neq \emptyset$ and, for some $T
\in \bigcap_{i = 1}^k \mathcal{T}_i$, $\mathcal{A}_i =
\mathcal{F}_i\langle T \rangle$ for all $i \in [k]$.\\
(ii) $\sum_{i = 1}^k |\mathcal{F}_i \langle T_i \rangle| <
\max\{|\mathcal{F}_1|, \dots, |\mathcal{F}_k|\}$
and, for some $j \in [k]$ such that $|\mathcal{F}_j| = \max
\{|\mathcal{F}_i| \colon i \in [k]\}$, $\mathcal{A}_j =
\mathcal{F}_j$ and $\mathcal{A}_i = \emptyset$ for all $i \in [k]
\backslash \{j\}$.\\
(iii) $\sum_{i = 1}^k |\mathcal{F}_i \langle T_i \rangle| =
\max\{|\mathcal{F}_1|, \dots, |\mathcal{F}_k|\}$,
and $\mathcal{A}_1, \dots, \mathcal{A}_k$ are as in (i) or (ii).
\end{theorem}
%Any conjectures similar to those in the previous subsection?

\begin{remark} \emph{Consider the special case where each of the sets $S_i$ in Theorem~\ref{sumresult} contains only one element $s_i$ (that
is, $\mathcal{F}_i = \mathcal{H}^{(s_i)}$ for all $i \in [k]$). Then,
by Lemma~\ref{Spernercor} below, for any $j \in [k]$, $|\mathcal{F}_j| = \max\{|\mathcal{F}_1|, \dots, |\mathcal{F}_k|\}$ if and only if $s_j = \max\{s_1, \dots, s_k\}$.}
\end{remark}

\noindent
\textbf{Proof of Theorem~\ref{sumcor2}.} Theorem~\ref{sumcor2} is the case $S_1 = \dots = S_k = S$, and hence $\mathcal{F}_1 = \dots = \mathcal{F}_k = \mathcal{F}$, in Theorem~\ref{sumresult}. Since $\mathcal{L}$ is a largest $t$-star of $\mathcal{F}$, the conditions of Theorem~\ref{sumresult} are satisfied.~\hfill{$\Box$}\\

One of the main challenges in establishing Theorem~\ref{sumresult} was to have a condition for $\mu(\mathcal{H})$ that is independent of $k$. Note that, if $k$ is sufficiently large, then part (i) of Theorem~\ref{sumresult} holds; in particular, this holds with $k \geq |\mathcal{H}^{(\leq r)}|$ because we then have $\sum_{i = 1}^k |\mathcal{F}_i \langle T_i \rangle| \geq \sum_{i = 1}^k 1 = k \geq |\mathcal{H}^{(\leq r)}| > \max\{|\mathcal{F}_1|, \dots, |\mathcal{F}_k|\}$. On the other hand, in Section~\ref{sumsection} we also prove the following result, which says that part (ii) holds if $\mu(\mathcal{H})$ is sufficiently large depending on $k$.

\begin{theorem} \label{sumcor1} If $\mathcal{H}$, $\mathcal{F}_1, \dots, \mathcal{F}_k$, $\mathcal{A}_1, \dots, \mathcal{A}_k$ are as in Theorem~\ref{sumresult}, $j \in [k]$ such that $|\mathcal{F}_i| \leq |\mathcal{F}_j|$ for all $i \in [k]$, and $\mu(\mathcal{H}) \geq \max\{n_{\rm S}(r,t), (k^{1/t}+1)r\}$, then
\[\sum_{i = 1}^k |\mathcal{A}_i| \leq |\mathcal{F}_j|,\]
and equality holds if and only if, for some $l \in [k]$ such that
$|\mathcal{F}_l| = |\mathcal{F}_j|$, $\mathcal{A}_l =
\mathcal{F}_l$ and $\mathcal{A}_i = \emptyset$ for all $i \in [k]
\backslash \{l\}$.
\end{theorem}

A \emph{graph} $G$ is a pair $(V,E)$ with $E \subseteq {V \choose 2}$, and a set $I \subseteq V$ is said to be an \emph{independent set of $G$} if $\{i,j\} \notin E$ for every $i, j \in I$. Let $\mathcal{I}_G$ denote the family of all independent sets of a graph $G$. Many EKR-type results can be phrased in terms of independent sets of graphs; see \cite[page 2878]{BH2}. Holroyd and Talbot \cite{HT} introduced the EKR problem for the families ${\mathcal{I}_G}^{(r)}$. They conjectured that ${\mathcal{I}_G}^{(r)}$ has the $1$-star property if $\mu(\mathcal{I}_G) \geq 2r$ \cite[Conjecture~7]{HT}. This conjecture inspired many results; see, for example, \cite{HST,HT,BH2,HK,Wr}. Clearly, $\mathcal{I}_G$ is a hereditary family, so Theorem~\ref{t int her} verifies the conjecture for $\mu(\mathcal{I}_G) \geq \frac{3}{2}(r-1)^2(3r-4)+r$. Kamat \cite{Kamat} made the following analogous conjecture for cross-intersecting families.

\begin{conj}[\cite{Kamat}] If $G$ is a graph, $\mu(\mathcal{I}_G) \geq 2r$, and $\mathcal{A}$ and $\mathcal{B}$ are cross-intersecting subfamilies of ${\mathcal{I}_G}^{(r)}$, then $|\mathcal{A}| + |\mathcal{B}| \leq |{\mathcal{I}_G}^{(r)}|$.
\end{conj}
We suggest the following generalization.

\begin{conj}\label{kamatgen} If $\mathcal{H}$ is a hereditary family, $\mu(\mathcal{H}) \geq 2r$, and $\mathcal{A}$ and $\mathcal{B}$ are cross-intersecting subfamilies of $\mathcal{H}^{(r)}$, then $|\mathcal{A}| + |\mathcal{B}| \leq |\mathcal{H}^{(r)}|$.
\end{conj}
We also conjecture that, if $\mu(\mathcal{H}) > 2r$, then the bound is attained if and only if one of $\mathcal{A}$ and $\mathcal{B}$ is $\mathcal{H}^{(r)}$ and the other is empty. This is true for $\mu(\mathcal{H})$ sufficiently large depending only on $r$.

\begin{theorem}\label{kamatproof} Conjecture~\ref{kamatgen} is true if $\mu(\mathcal{H}) \geq n_{\rm S}(r,1)$.
\end{theorem}
\textbf{Proof.} This is Theorem~\ref{sumcor1} with $t=1$, $k = 2$ and $S_1 = S_2 = \{r\}$, in which case $\max\{n_{\rm S}(r,t), (k^{1/t}+1)r\} = n_{\rm S}(r,t)$.~\hfill{$\Box$} \\

For the special case where $\mathcal{H}$ is the power set of $[n]$, Theorem~\ref{sumresult} yields the following.

\begin{theorem} \label{sumcor3} Let $1 \leq t \leq r$, $n \geq n_{\rm S}(r,t)$ and $k \geq 2$. For all $i \in [k]$, let $\emptyset \neq S_i \subseteq [t,r]$ and $\mathcal{F}_i := \bigcup_{s \in S_i}{[n]
\choose s}$. Let $\mathcal{A}_1, \dots, \mathcal{A}_k$ be cross-$t$-intersecting families such that $A_i \subseteq \mathcal{F}_i$ for all $i \in [k]$. Then
\[ \sum_{i = 1}^k |\mathcal{A}_i| \leq \max \left\{ \sum_{i = 1}^k
\sum_{s \in S_i}{n-t \choose s-t}, \sum_{s \in S_1}{n \choose s},
\dots, \sum_{s \in S_k}{n \choose s} \right\}.\]
Moreover, equality holds if and only if one of the following
holds: \\
(i) $\sum_{i = 1}^k \sum_{s \in S_i}{n-t \choose s-t} > \max
\left\{\sum_{s \in S_i}{n \choose s} \colon i \in [k] \right\}$
and, for some $t$-subset $T$ of $[n]$, $\mathcal{A}_i =
\mathcal{F}_i \langle T \rangle$ for all $i \in [k]$.\\
(ii) $\sum_{i = 1}^k \sum_{s \in S_i}{n-t \choose s-t} <
\max\left\{\sum_{s \in S_i}{n \choose s} \colon i \in [k]
\right\}$ and, for some $j \in [k]$ such that $|\mathcal{F}_j| =
\max \{|\mathcal{F}_i| \colon i \in [k]\}$, $\mathcal{A}_j =
\mathcal{F}_j$ and $\mathcal{A}_i = \emptyset$ for all $i \in [k]
\backslash \{j\}$.\\
(iii) $\sum_{i = 1}^k \sum_{s \in S_i}{n-t \choose s-t} =
\max\left\{\sum_{s \in S_i}{n \choose s} \colon i \in [k]
\right\}$, and $\mathcal{A}_1, \dots, \mathcal{A}_k$ are as in (i) or~(ii).
\end{theorem}
\textbf{Proof.} For any $t$-subset $T$ of $[n]$, $\mathcal{F}_i
\langle T \rangle$ is a largest $t$-star of $\mathcal{F}_i$,
because all the $t$-stars of $\mathcal{F}_i$ are of size $\sum_{s
\in S_i}{n-t \choose s-t}$. Let $\mathcal{H} = 2^{[n]}$. Then $n =
\mu(\mathcal{H})$ and ${[n] \choose s} = \mathcal{H}^{(s)}$ (for
any $s$). The result now follows immediately from
Theorem~\ref{sumresult}.~\hfill{$\Box$} \\

As pointed out in Section~\ref{Intro}, the case where $S_1 = \dots = S_k = \{r\}$, and hence $\mathcal{F}_1 = \dots = \mathcal{F}_k = {[n] \choose r}$, was settled for every $n$ by Wang and Zhang \cite{WZ}, and the case $t=1$ of their result is Hilton's seminal result \cite{H}.

\subsection{The maximum product} \label{prodsubsec}

The following is the product version of Theorem~\ref{sumresult}. It was proved in \cite{Borg_maxprod} and is proved in a different way in Section~\ref{prodsection}. It is our main and most general result for the maximum product.

\begin{theorem} \label{prodresult} Let $1 \leq t \leq r$, and let $\mathcal{H}$ be a hereditary family with $\mu(\mathcal{H}) \geq n_{\rm P}(r,t)$. Let $k \geq 2$, and, for each $i \in [k]$, let $\emptyset \neq S_i \subseteq [t,r]$, $\mathcal{F}_i := \bigcup_{s \in
S_i}\mathcal{H}^{(s)}$ and $\mathcal{T}_i := \{T \in \mathcal{H}^{(t)} \colon \mathcal{F}_i \langle T \rangle$ is a largest
$t$-star of $\mathcal{F}_i\}$. Let $\mathcal{A}_1, \dots, \mathcal{A}_k$ be cross-$t$-intersecting families such that $\mathcal{A}_i
\subseteq \mathcal{F}_i$ for all $i \in [k]$. Then $\mathcal{T}_i
\neq \emptyset$ for all $i \in [k]$, and, for any $T_1 \in
\mathcal{T}_1, \dots, T_k \in \mathcal{T}_k$,
$$\prod_{i = 1}^k |\mathcal{A}_i| \leq
\prod_{i = 1}^k |\mathcal{F}_i \langle T_i \rangle|.$$
Moreover, equality holds if and only if $\bigcap_{i = 1}^k
\mathcal{T}_i \neq \emptyset$ and, for some $T \in \bigcap_{i =
1}^k \mathcal{T}_i$, $\mathcal{A}_i = \mathcal{F}_i\langle T
\rangle$ for all $i \in [k]$.
\end{theorem}
We conjecture that there exists an integer $n^*_{\rm P}(r,t)$ such that, if $\mu(\mathcal{H}) \geq n^*_{\rm P}(r,t)$, then $\prod_{i = 1}^k |\mathcal{A}_i| \leq \prod_{i = 1}^k |\mathcal{F}_i \langle T \rangle|$ for some $T \in \mathcal{H}^{(t)}$. Clearly, by Theorem~\ref{prodresult}, this is true with $n^*_{\rm P}(r,t) = n_{\rm P}(r,t)$ if $\bigcap_{i = 1}^k \mathcal{T}_i \neq \emptyset$; however, \cite[Example 1]{Borg} shows that $\bigcap_{i = 1}^k \mathcal{T}_i$ may be empty. In view of Theorem~\ref{AK}, we also conjecture that the conjecture holds with $n^*_{\rm P}(r,t) = (t+1)(r-t+1)$.

The remaining results are analogues of results in Section~\ref{maxsumsubsec}. Next, we have the result for cross-$t$-intersecting subfamilies of a union of levels of a hereditary family.

\begin{theorem} \label{prodcor1} If $1 \leq t \leq
r$, $\mathcal{H}$ is a hereditary family with $\mu(\mathcal{H}) \geq n_{\rm P}(r,t)$, $\emptyset \neq S \subseteq [t,r]$, $\mathcal{F} := \bigcup_{s \in S}\mathcal{H}^{(s)}$, $\mathcal{A}_1, \dots, \mathcal{A}_k$ are cross-$t$-intersecting subfamilies of
$\mathcal{F}$, and $\mathcal{L}$ is a largest $t$-star of
$\mathcal{F}$, then
\[\prod_{i = 1}^k |\mathcal{A}_i| \leq |\mathcal{L}|^k,\]
and equality holds if and only if $\mathcal{A}_1 = \dots = \mathcal{A}_k = \mathcal{L}'$ for some largest $t$-star $\mathcal{L}'$ of $\mathcal{F}$.
\end{theorem}
\textbf{Proof.} This is the case $S_1 = \dots = S_k = S$, and
hence $\mathcal{F}_1 = \dots = \mathcal{F}_k =
\mathcal{F}$, in Theorem~\ref{prodresult}. Since $\mathcal{L}$ is
a largest $t$-star of $\mathcal{F}$, the conditions of
Theorem~\ref{prodresult} are satisfied.~\hfill{$\Box$} \\

Similarly to Theorem~\ref{sumcor2}, the result above implies Theorem~\ref{AK gen} for $\mu(\mathcal{H}) \geq n_{\rm P}(r,t)$. Indeed, if $\mathcal{A}$ is an intersecting subfamily of $\mathcal{F}$ and $\mathcal{A}_1 = \dots = \mathcal{A}_k = \mathcal{A}$, then $\mathcal{A}_1, \dots, \mathcal{A}_k$ are cross-$t$-intersecting subfamilies of $\mathcal{F}$, so $|\mathcal{A}| = \left(\prod_{i = 1}^k |\mathcal{A}_i|\right)^{1/k} \leq \left(|\mathcal{L}|^k \right)^{1/k} = |\mathcal{L}|$.

\begin{theorem} \label{prodcor1.5} If $1 \leq t \leq r$, $\mathcal{H}$ is a hereditary family with $\mu(\mathcal{H}) \geq n_{\rm P}(r,t)$, $\mathcal{A}_1, \dots, \mathcal{A}_k$ are cross-$t$-intersecting subfamilies of $\mathcal{H}^{(\leq r)}$, and $\mathcal{L}$ is a largest $t$-star of $\mathcal{H}^{(\leq r)}$, then
\[\prod_{i = 1}^k |\mathcal{A}_i| \leq |\mathcal{L}|^k,\]
and equality holds if and only if $\mathcal{A}_1 = \dots = \mathcal{A}_k = \mathcal{L}'$ for some largest $t$-star $\mathcal{L}'$ of $\mathcal{H}^{(\leq r)}$.
\end{theorem}
\textbf{Proof.} If one of the families $\mathcal{A}_1, \dots, \mathcal{A}_k$ has a member of size less than $t$, then, since $\mathcal{A}_1, \dots, \mathcal{A}_k$ are cross-$t$-intersecting, the other families are empty, and hence $\prod_{i = 1}^k |\mathcal{A}_i| = 0$. If none of $\mathcal{A}_1, \dots, \mathcal{A}_k$ has a member of size less than $t$, then $\mathcal{A}_1, \dots, \mathcal{A}_k \subseteq \bigcup_{s \in [t,r]}\mathcal{H}^{(s)}$, and hence the result follows by Theorem~\ref{prodcor1}.~\hfill{$\Box$} \\

The problem for $\mathcal{H}^{(\leq r)}$ for the special case where $\mathcal{H} = 2^{[n]}$ and $t=1$ is solved in \cite{Borg6} for every $r$.

We conclude this section with an analogue of Theorem~\ref{sumcor3}.

\begin{theorem} \label{prodcor2} Let $1 \leq t \leq
r$, $n \geq n_{\rm P}(r,t)$ and $k \geq 2$. For all $i \in [k]$, let $\emptyset \neq S_i \subseteq [t,r]$ and $\mathcal{F}_i := \bigcup_{s \in S_i}{[n] \choose s}$. If $\mathcal{A}_1, \dots,
\mathcal{A}_k$ are cross-$t$-intersecting families such that $A_i
\subseteq \mathcal{F}_i$ for all $i \in [k]$, then
$$\prod_{i = 1}^k |\mathcal{A}_i| \leq
\prod_{i = 1}^k |\mathcal{F}_i \langle [t] \rangle| = \prod_{i =
1}^k \sum_{s \in S_i}{n-t \choose s-t},$$
and equality holds if and only if, for some $T \in {[n] \choose t}$, $\mathcal{A}_i = \mathcal{F}_i\langle T \rangle$ for all
$i \in [k]$.
\end{theorem}
\textbf{Proof.} For each $T \in {[n] \choose t}$, $\mathcal{F}_i
\langle T \rangle$ is a largest $t$-star of $\mathcal{F}_i$ as all the $t$-stars of $\mathcal{F}_i$ are of size $\sum_{s
\in S_i}{n-t \choose s-t}$. Let $\mathcal{H} = 2^{[n]}$. Then $n =
\mu(\mathcal{H})$ and ${[n] \choose s} = \mathcal{H}^{(s)}$ (for
any $s$). The result now follows immediately from
Theorem~\ref{prodresult}.~\hfill{$\Box$} \\

%Theorem~\ref{prodcor2} is proved in \cite{Borg11} also for $n$ sufficiently large. 
Theorem~\ref{prodcor2} is proved in \cite{Borg12} with a condition on $n$ that is close to best possible. In each of \cite{Tok1,Tok2,FLST}, Theorem~\ref{prodcor2} is proved for $k = 2$ and $S_1 = S_2 = \{r\}$ with a condition on $n$ that is also nearly optimal.
%close to best possible.

\section{An intersection lemma} \label{keysection}

We now start working towards the proofs of Theorems~\ref{sumresult}, \ref{sumcor1} and \ref{prodresult}. In this section, we focus on the structure of families that are not trivial $t$-intersecting families, and on $t$-intersections with their sets. We start by establishing the following key lemma, which is best possible.

\begin{lemma}\label{keyint} If $r$ is the size of a largest set in a family $\mathcal{A}$, $|\mathcal{A}| > 1$, and $\mathcal{A}$ is not a trivial $t$-intersecting family, then, for some integer $p$ with $2 \leq p \leq \max\{2,r-t+2\}$, $\mathcal{A}$ has $p$ sets that have a union of size at most $m(r,t)$ and do not have $t$ common elements.
\end{lemma}
\textbf{Proof.} If $\mathcal{A}$ is not $t$-intersecting, then there exist two sets $A$ and $B$ in $\mathcal{A}$ such that $|A \cap B| < t$, and hence the result is immediate since $|A \cup B| \leq |A| + |B|
\leq 2r$.

Now suppose that $\mathcal{A}$ is a $t$-intersecting family but not a trivial one. Let $A_1$ and $A_2$ be two distinct sets in $\mathcal{A}$. Then $|A_1 \cap A_2| \geq t$. Since $\mathcal{A}$ is not trivial, we can find a set $A_3$ in $\mathcal{A}$ such that $A_1 \cap A_2 \nsubseteq A_3$, so
$|A_1 \cap A_2 \cap A_3| < |A_1 \cap A_2|$. If $|A_1 \cap A_2
\cap A_3| \geq t$, then, since $\mathcal{A}$ is not trivial, we
can find a set $A_4$ in $\mathcal{A}$ such that $A_1 \cap A_2
\cap A_3 \nsubseteq A_4$, so $|A_1 \cap A_2 \cap A_3 \cap
A_4| < |A_1 \cap A_2 \cap A_3|$. Continuing this way, we eventually
obtain $p$ sets $A_1, A_2, \dots, A_p$ in $\mathcal{A}$ such that
\begin{gather} |A_1 \cap A_2 \cap \dots \cap A_p| < |A_1 \cap A_2 \cap \dots \cap A_{p-1}| < \dots < |A_1|, \label{key1} \\
|A_1 \cap A_2 \cap \dots \cap A_p| < t \leq |A_1 \cap A_2 \cap \dots
\cap A_{p-1}|, \label{key2}
\end{gather}
where $3 \leq p \leq r-t+2$.

We next show that $|A_1 \cup A_{p-1} \cup A_p| \leq 3r - 2t - 1$.
Since $\mathcal{A}$ is $t$-intersecting, we have
$|A_p \cap A_1| \geq t$, $|A_{p-1} \cap A_1| \geq t$ and $|A_p
\cap A_{p-1}| \geq t$.

Suppose $|A_{p-1} \cap A_1| = t$. Then, by (\ref{key2}), $A_1
\cap A_2 \cap \dots \cap A_{p-1} = A_{p-1} \cap A_1$, and hence $|A_p
\cap (A_{p-1} \cap A_1)| = |A_1 \cap A_2 \cap \dots \cap A_p| \leq
t-1$. Thus, since $t \leq |A_p \cap A_{p-1}| = |A_p \cap (A_{p-1}
\cap A_1)| + |A_p \cap (A_{p-1} \backslash A_1)|$, we get $|A_p
\cap (A_{p-1} \backslash A_1)| \geq 1$. Thus, $A_p$ intersects $A_1
\cup A_{p-1}$ in at least $t+1$ elements, and hence $|A_p
\backslash (A_1 \cup A_{p-1})| \leq r - t - 1$. We have $|A_1
\cup A_{p-1} \cup A_p| = |A_p \backslash (A_1 \cup A_{p-1})| +
|A_1 \cup A_{p-1}| = |A_p \backslash (A_1 \cup A_{p-1})| + (|A_1|
+ |A_{p-1}| - |A_1 \cap A_{p-1}|) \leq (r-t-1) + (2r - t) = 3r -
2t - 1$.

Now suppose $|A_{p-1} \cap A_1| \geq t+1$. Then $|A_{p-1}
\backslash (A_1 \cup A_p)| \leq r - t - 1$. We have $|A_1 \cup
A_{p-1} \cup A_p| = |A_{p-1} \backslash (A_1 \cup A_p)| + |A_1
\cup A_p| \leq (r-t-1) + (2r - t) = 3r - 2t - 1$.

Recall that $p \geq 3$. Suppose $p = 3$. Then $A_1 \cup
A_{p-1} \cup A_p = A_1 \cup A_2 \cup A_3$. Thus, $|A_1 \cup A_2 \cup
A_3| \leq 3r-2t-1$. Since $|A_1 \cap A_2| \geq t$ and $A_1 \neq A_2$, $r \geq t+1$. Thus, since $\left(\frac{1}{2}(r-t)(r-t+5) + t-1 \right) - (3r-2t-1) = \frac{1}{2}(r-t)(r-t-1)$, we have $3r-2t-1 \leq
\frac{1}{2}(r-t)(r-t+5) + t-1$, as required.

Now suppose $p \geq 4$. Consider any set $A_q$ with $2 \leq q \leq
p-2$. By (\ref{key1}) and (\ref{key2}), $t \leq |A_1 \cap A_2 \cap \dots \cap A_{p-1}|$, $t+1 \leq |A_1 \cap A_2 \cap \dots \cap
A_{p-2}|$, $\dots$, $t + p - 1 - q \leq |A_1 \cap A_2 \cap \dots \cap A_q|$. Since $|A_1 \cap A_q| \geq |A_1 \cap A_2 \cap \dots \cap A_q| \geq t+p-1-q$, $|A_q \backslash A_1| \leq r - (t+p-1-q) = r - t - p + 1 + q$. We have

\begin{align} |A_1 \cup A_2 \cup \dots \cup A_p| & \leq |A_1 \cup
A_{p-1} \cup A_p| + |(A_2 \cup \dots \cup A_{p-2}) \backslash (A_1
\cup A_{p-1} \cup A_p)| \nonumber \\
& \leq |A_1 \cup A_{p-1} \cup A_p| + |(A_2 \cup \dots \cup A_{p-2})
\backslash A_1| \nonumber \\
&\leq |A_1 \cup A_{p-1} \cup A_p| + |A_2 \backslash A_1| + \dots +
|A_{p-2} \backslash A_1| \nonumber \\
&\leq (3r - 2t - 1) + (r - t - p + 3) + \dots + (r - t - 1)
\nonumber \\
&= (3r - 2t - 1) + \sum_{i = 1}^{p-3} (r-t-i) \nonumber \\
&\leq (3r - 2t - 1) + \sum_{i = 1}^{r-t-1} (r-t-i) \quad
\mbox{(since $p \leq r-t+2$)} \nonumber \\
&= (3r - 2t - 1) + \sum_{j = 1}^{r-t-1} j = 2(r-t) + r - 1 +
\frac{(r-t-1)(r-t)}{2} \nonumber \\
%&= \frac{4(r-t) + 2r - 2 + (r-t)^2 - (r-t)}{2} = \frac{(r-t)^2 +
%5(r-t) + 2t - 2}{2} \nonumber \\
&= %\frac{(r-t)^2 + 5(r-t)}{2} + t-1
\frac{(r-t)(r-t+5)}{2} +
t-1, \nonumber
\end{align}
as required.~\hfill{$\Box$}

\begin{remark} \emph{The upper bound on the size of the union in the
lemma above is sharp. An example for the case $2r \geq
\frac{1}{2}(r-t)(r-t+5) + t-1$ is when $\mathcal{A}$ consists of
two disjoint $r$-sets. Now consider $2r \leq
\frac{1}{2}(r-t)(r-t+5) + t-1$. Then $r \geq t+2$. Let $x_{0,1}, \dots, x_{0,r}$, $x_{1,1}$, $x_{2,1}, x_{2,2}, \dots, x_{r-t,1}, \dots,
x_{r-t,r-t}$, $y_1, \dots, y_{r-t-1}$ be distinct numbers. Let
$A_0 := \{x_{0,j} \colon j \in [r]\}$. For each $i \in [r-t]$, let
$A_i := \{x_{0,j} \colon j \in [r-i]\} \cup \{x_{i,j} \colon j \in
[i]\}$. Let $A_{r-t+1} := \{x_{0,j} \colon j \in [t-1]\} \cup \{x_{0,t+1}, x_{r-t,1}\} \cup \{y_j \colon j \in [r-t-1]\}$.
Now let $\mathcal{A} := \{A_i \colon i \in \{0\} \cup [r-t+1]\}$. It is easy to check that $|\bigcap_{i=0}^{r-t+1} A_i| = t-1$ and
$|\bigcup_{i=0}^{r-t+1} A_i| = \frac{1}{2}(r-t)(r-t+5) + t-1$.}
\end{remark}

A set that $t$-intersects each set in a family $\mathcal{A}$ is called a \emph{$t$-transversal of $\mathcal{A}$}.

\begin{lemma}\label{key1corA} If $r$ is the size of a largest set in a family $\mathcal{A}$, $|\mathcal{A}| > 1$, and $\mathcal{A}$ is not a trivial $t$-intersecting family, then, for some $p \geq 2$, there exist $A_1, \dots, A_p \in \mathcal{A}$ such that $|\bigcup_{i=1}^p A_i| \leq m(r,t)$ and every $t$-transversal of $\mathcal{A}$ $(t+1)$-intersects $\bigcup_{i=1}^p A_i$.
\end{lemma}
\textbf{Proof.} Lemma~\ref{keyint} tells us that, under the given conditions, there
exist $p \geq 2$ sets $A_1, \dots, A_p$ in $\mathcal{A}$ such that
$|\bigcap_{i=1}^p A_i| \leq t-1$ and $|\bigcup_{i=1}^p A_i| \leq m(r,t)$. Let $C := \bigcup_{i=1}^p A_i$. Let $B$ be a $t$-transversal of $\mathcal{A}$. Let $D := B \cap C$. Then $|D| \geq |B \cap A_1| \geq t$. Suppose $|D| = t$. Since $|\bigcap_{i=1}^p A_i| \leq t-1$, $D \nsubseteq A_j$ for some $j \in [p]$. We obtain $|D \cap A_j| \leq t-1$, which contradicts $|D \cap A_j| = |(B \cap C) \cap A_j| = |B \cap (C \cap A_j)| = |B \cap A_j|  \geq t$. Thus, $|D| \geq t+1$.~\hfill{$\Box$}

%\section{Key properties of hereditary families}
\section{Results on hereditary families}
\label{propertysection}

The results provided in this section establish the properties of hereditary families that are needed for the proofs of our main results. The first one is given by {\cite[Corollary~3.2]{Borg}}. 

\begin{lemma}[{\cite{Borg}}]\label{Spernercor} If $\mathcal{H}$ is a hereditary family and $0 \leq p < q \leq \mu(\mathcal{H}) - p$, then
\[|\mathcal{H}^{(q)}| \geq \frac{{\mu(\mathcal{H}) - p \choose
q - p}}{{q \choose q-p}}|\mathcal{H}^{(p)}|.\]
\end{lemma}

We will need the next lemma (for arbitrary families) in the result that follows it.

\begin{lemma} \label{mucor} If $\mathcal{F}$ is a family and $X$ is a set such that $\mathcal{F} \langle X \rangle \neq \emptyset$, then
\[\mu(\{F \backslash X \colon F \in \mathcal{F} \langle X \rangle\}) \geq \mu(\mathcal{F}) - |X|.\]
\end{lemma}
\textbf{Proof.} Let $\mathcal{G} := \{F \backslash X \colon F \in \mathcal{F} \langle X \rangle\}$. Let $B$ be a base of $\mathcal{G}$ of size $\mu(\mathcal{G})$. Then $B \cup X$ is a base of $\mathcal{F}$. Thus, $\mu(\mathcal{F}) \leq |B| + |X| = \mu(\mathcal{G}) + |X|$.~\hfill{$\Box$}

\begin{lemma}\label{mainlemma2} If $0 \leq t_1 < t_2 \leq r$, $\mathcal{H}$ is a hereditary family with $\mu(\mathcal{H}) \geq 2r - t_1$, $\emptyset \neq S \subseteq [t_2,r]$, $\mathcal{F} := \bigcup_{s \in S}\mathcal{H}^{(s)}$, and $T_1$ is a $t_1$-subset of a $t_2$-set $T_2$ such that $\mathcal{F}\langle T_2 \rangle \neq \emptyset$, then
\[|\mathcal{F}\langle T_1 \rangle| > \frac{{\mu(\mathcal{H}) - r \choose t_2 - t_1}} {{r-t_1 \choose t_2 - t_1}} |\mathcal{F} \langle
T_2 \rangle|.\]
\end{lemma}
\textbf{Proof.} Let $s \in S$. Let $F \in \mathcal{F}\langle T_2
\rangle$, and let $G$ be a base of $\mathcal{H}$ such that $F \subseteq G$. Since $|G| \geq \mu(\mathcal{H}) \geq r \geq s$ and $\mathcal{H}$ is hereditary, $\emptyset \neq {G
\choose s} \subseteq \mathcal{H}^{(s)}$. Since $T_2 \subseteq F
\subseteq G$ and $|T_2| = t_2 \leq s$, we then have
$\mathcal{H}^{(s)}\langle T_2 \rangle \neq \emptyset$. Thus,
$\mathcal{H} \langle T_2 \rangle \neq \emptyset$.

Let $\mathcal{I} := \{H \backslash T_2 \colon H \in
\mathcal{H}\langle T_2 \rangle\}$. Since $\mathcal{H}$ is
hereditary, $\mathcal{I}$ is hereditary. By Lemma~\ref{mucor},
\begin{equation} \mu(\mathcal{I}) \geq \mu(\mathcal{H}) - t_2.
\label{t int her 2}
\end{equation}
Let $p := s - |T_2| = s - t_2$, $q := p + t_2 - t_1 = s-t_1$.
Given that $\mu(\mathcal{H}) \geq 2r - t_1$, it follows by (\ref{t
int her 2}) that
\begin{equation} \mu(\mathcal{I}) \geq 2r - t_1 - t_2 \geq 2s -
t_1 - t_2 = p + q. \nonumber
\end{equation}
Therefore, by Lemma~\ref{Spernercor}, we have
\begin{gather} |\mathcal{I}^{(q)}| \geq
\frac{ {\mu(\mathcal{I})-p \choose t_2 - t_1} }{ {q \choose t_2 -
t_1} } |\mathcal{I}^{(p)}|, \nonumber
\end{gather}
and hence, since $|\mathcal{I}^{(p)}| = |\mathcal{H}^{(s)}\langle
T_2 \rangle|$ and $|\mathcal{I}^{(q)}| = |\mathcal{H}^{(s + t_2 -
t_1)} \langle T_2 \rangle|$ (by definition of $\mathcal{I}$, $p$
and $q$),
\begin{align} |\mathcal{H}^{(s + t_2 - t_1)} \langle T_2 \rangle|
&\geq \frac{ {\mu(\mathcal{I})-p \choose t_2 - t_1} }{ {q \choose
t_2 - t_1} } |\mathcal{H}^{(s)}\langle T_2 \rangle| \nonumber \\
&\geq \frac{ {(\mu(\mathcal{H}) -t_2) - (s - t_2) \choose t_2 -
t_1} }{ {s-t_1 \choose t_2 - t_1} } |\mathcal{H}^{(s)}\langle T_2
\rangle| \quad \mbox{(by (\ref{t int her 2}))} \nonumber \\
&= \frac{ {\mu(\mathcal{H})-s \choose t_1 - t_2} }{ {s-t_1 \choose
t_2 - t_1} } |\mathcal{H}^{(s)}\langle T_2 \rangle| \geq \frac{
{\mu(\mathcal{H})- r \choose t_2 - t_1} }{ {r-t_1 \choose t_2 - t_1} }
|\mathcal{H}^{(s)}\langle T_2 \rangle|. \label{t int her 3}
\end{align}
Let $T_3 := T_2 \backslash T_1$ and $\mathcal{B} := \{A
\backslash T_3 \colon A \in \mathcal{H}^{(s + t_2 - t_1)}\langle
T_2 \rangle\}$. For every $B \in \mathcal{B}$, we have $T_1 \subset
B$ and $|B| = s$. Thus, since $\mathcal{H}$ is hereditary,
$\mathcal{B} \subseteq \mathcal{H}^{(s)}\langle T_1 \rangle$.
Since $\emptyset \neq \mathcal{H}^{(s)} \langle T_2 \rangle
\subseteq \mathcal{H}^{(s)} \langle T_1 \rangle \backslash
\mathcal{B}$, we actually have $\mathcal{B} \subsetneq
\mathcal{H}^{(s)}\langle T_1 \rangle$, and hence $|\mathcal{B}| <
|\mathcal{H}^{(s)} \langle T_1 \rangle|$. Thus, since
$|\mathcal{B}| = |\mathcal{H}^{(s + t_2 - t_1)} \langle T_2
\rangle|$, $|\mathcal{H}^{(s + t_2 - t_1)} \langle T_2
\rangle| < |\mathcal{H}^{(s)} \langle T_1 \rangle|$. Therefore, by
(\ref{t int her 3}),
\begin{equation} |\mathcal{H}^{(s)} \langle T_1 \rangle| >
\frac{ {\mu(\mathcal{H})- r \choose t_2 - t_1} }{ {r-t_1 \choose t_2
- t_1} } |\mathcal{H}^{(s)}\langle T_2 \rangle|. \label{t int her
4}
\end{equation}

The result follows immediately from (\ref{t int her 4}) since
$|\mathcal{F} \langle T_1 \rangle| = \sum_{s \in S}
|\mathcal{H}^{(s)} \langle T_1 \rangle|$ and $|\mathcal{F} \langle
T_2 \rangle| = \sum_{s \in S} |\mathcal{H}^{(s)} \langle T_2
\rangle|$.~\hfill{$\Box$}
\\

%The new result above has the following two consequences, the first of which is a new result that is needed for the sum results, and the second of which is a result in \cite{Borg} that is needed for both the sum and the product results.
The new result above has the following two consequences, the first of which is needed for the sum results, and the second of which is needed for both the sum and the product results.

\begin{lemma}\label{maincor2} If $1 \leq t \leq r$, $\mathcal{H}$ is a hereditary family with $\mu(\mathcal{H}) \geq 2r$,  $\emptyset \neq S \subseteq [t,r]$, $\mathcal{F} := \bigcup_{s \in S}\mathcal{H}^{(s)}$, and $\mathcal{L}$ is a largest $t$-star of $\mathcal{F}$, then
\[|\mathcal{F}| > \frac{{\mu(\mathcal{H}) - r \choose t}}{{r
\choose t}} |\mathcal{L}|.\]
\end{lemma}
\textbf{Proof.} Since $\mathcal{H}$ is hereditary and
$\mu(\mathcal{H}) \geq 2r > t$, we obviously have $\mathcal{F}
\neq \emptyset$, and hence $\mathcal{L} \neq \emptyset$. Thus,
$\mathcal{L} = \mathcal{F} \langle T \rangle \neq \emptyset$ for
some $t$-set $T$. Since $\mathcal{F} = \mathcal{F} \langle
\emptyset \rangle$, the result now follows by
Lemma~\ref{mainlemma2} with $T_1 = \emptyset$ and $T_2 = T$.~\hfill{$\Box$}

\begin{lemma} %[\cite{Borg}] 
\label{main lemma} If $1 \leq t+1 \leq r$, $\emptyset \neq S \subseteq [t+1,r]$, $\mathcal{H}$ is a hereditary
family with $\mu(\mathcal{H}) \geq 2r-t$, $\mathcal{F} :=
\bigcup_{s \in S}\mathcal{H}^{(s)}$, $\emptyset \neq \mathcal{A} \subseteq \mathcal{F}$, and $X$ is a $(t+1)$-transversal of $\mathcal{A}$, then there exists some $T \in {X \choose t}$ such that
\[|\mathcal{A}| < \frac{r-t}{\mu(\mathcal{H}) - r} {|X| \choose t+1} |\mathcal{F} \langle T \rangle |.\]
\end{lemma}
\textbf{Proof.} Choose $I_0 \in {X \choose t+1}$ such that
$|\mathcal{F}\langle I \rangle| \leq |\mathcal{F} \langle I_0
\rangle|$ for all $I \in {X \choose t+1}$. Given that $\emptyset
\neq \mathcal{A} \subseteq \mathcal{F}$ and $|A \cap X| \geq t+1$
for all $A \in \mathcal{A}$, we have
\begin{align} 1 \leq |\mathcal{A}| &= \left| \bigcup_{I \in {X
\choose t+1}}\mathcal{A}\langle I \rangle \right| \leq \sum_{I \in
{X \choose t+1}}|\mathcal{A}\langle I \rangle| \leq \sum_{I \in {X
\choose t+1}}|\mathcal{F} \langle I \rangle| \nonumber \\
&\leq \sum_{I \in {X \choose t+1}} |\mathcal{F} \langle I_0 \rangle| = {|X| \choose t+1} |\mathcal{F} \langle I_0 \rangle|. \nonumber
\end{align}
Choose $i_0 \in I_0$, and let $T := I_0 \backslash \{i_0\}$. By
Lemma~\ref{mainlemma2} with $T_1 = T$ and $T_2 = I_0$, we have
$|\mathcal{F} \langle T \rangle| > \frac{\mu(\mathcal{H}) -
r}{r-t} |\mathcal{F} \langle I_0 \rangle|$, so
$|\mathcal{F} \langle I_0 \rangle| < \frac{r-t}{\mu(\mathcal{H}) -
r} |\mathcal{F} \langle T \rangle|$. The result follows.~\hfill{$\Box$}

\section{Proofs of Theorems~\ref{sumresult} and \ref{sumcor1}} \label{sumsection}

We can now prove Theorems~\ref{sumresult} and \ref{sumcor1}, which yielded all the other results in Section~\ref{maxsumsubsec}. %our main result for the sum.\\
\\
\\
\textbf{Proof of Theorem~\ref{sumresult}.} Since $\mathcal{H}$ is hereditary and $\mu(\mathcal{H}) \geq n_{\rm S}(r,t) \geq r$, we have $\mathcal{F}_i \neq \emptyset$ and $\mathcal{T}_i \neq \emptyset$ for all $i \in [k]$. For each $i \in [k]$, let $T_i \in \mathcal{T}_i$. It is obvious that $\sum_{i = 1}^k |\mathcal{A}_i| = \max\left\{\sum_{i = 1}^k |\mathcal{F}_i \langle T_i \rangle|, |\mathcal{F}_1|, \dots, |\mathcal{F}_k|\right\}$ if one of (i)--(iii) holds. We now show that $\sum_{i = 1}^k
|\mathcal{A}_i| \leq \max\left\{\sum_{i = 1}^k |\mathcal{F}_i
\langle T_i \rangle|, |\mathcal{F}_1|, \dots, |\mathcal{F}_k|\right\}$, and that equality holds only if one of (i)--(iii) holds.

Suppose $r = t$. Then every member of $\mathcal{A}_i$ is
a $t$-set. Suppose that, for some $j \in [k]$, $\mathcal{A}_j$ is
non-empty. If $\mathcal{A}_j$ has two distinct members $A$ and
$B$, then no $t$-set $C$ can satisfy $|A \cap C| \geq t$ and $|B
\cap C| \geq t$, and hence, by the cross-$t$-intersection
condition, $\mathcal{A}_i = \emptyset$ for each $i \in [k]
\backslash \{j\}$. If $\mathcal{A}_j$ has only one member $A$,
then, by the cross-$t$-intersection condition, $\mathcal{A}_i \subseteq \{A\}$ for each $i \in [k] \backslash \{j\}$. Thus, the result for this case is trivial.

We now consider the case $r \geq t+1$. We shall abbreviate
$\mu(\mathcal{H})$, $m(r,t)$ and $n_{\rm S}(r,t)$ to $\mu$, $m$
and $n_{\rm S}$, respectively.\medskip

\textit{Case 1: For some $j \in [k]$, $\mathcal{A}_j$ is
non-empty and not a trivial $t$-intersecting family.}  Then $|\mathcal{A}_j| \geq 2$ because $|A| \geq t$ for all $A \in \mathcal{A}_j$ (as $\mathcal{A}_j \subseteq \mathcal{F}_j$). Let $h \in [k] \backslash \{j\}$. By Lemma~\ref{key1corA} and the cross-$t$-intersection
condition, there exists some $C \subseteq \bigcup_{A \in \mathcal{A}_j} A$ such that $|C| \leq m$ and $|B \cap C| \geq
t+1$ for all $B \in \mathcal{A}_h$. By Lemma~\ref{main lemma},
there exists some $T_0 \in {C \choose t}$ such that
\[|\mathcal{A}_h| < \frac{r-t}{\mu - r} {|C| \choose
t+1} |\mathcal{F}_h \langle T_0 \rangle |.\]
Since $|C| \leq m$ and
$\mathcal{F}_h\langle T_0 \rangle$ is a $t$-star of
$\mathcal{F}_h$, we then have
\[|\mathcal{A}_h| < \frac{r-t}{\mu - r} {m \choose
t+1} |\mathcal{F}_h\langle T_h \rangle|.\]
If $\mathcal{A}_h$ has a member $D$, then, since $|A \cap D| \geq t$ for all $A \in \mathcal{A}_j$,
\begin{align} |\mathcal{A}_j| &= \left|\bigcup_{E \in {D \choose
t}}\mathcal{A}_j \langle E \rangle \right| \leq \sum_{E \in {D
\choose t}}|\mathcal{A}_j \langle E \rangle| \leq \sum_{E \in {D
\choose t}} |\mathcal{F}_j \langle E \rangle| \leq \sum_{E \in {D
\choose t}}|\mathcal{F}_j \langle T_j \rangle| \nonumber \\
&= {|D| \choose t}|\mathcal{F}_j \langle T_j \rangle| \leq {r
\choose t}|\mathcal{F}_j \langle T_j \rangle|. \nonumber
\end{align}

We have shown that
\[|\mathcal{A}_i| < \frac{r-t}{\mu - r} {m \choose t+1}
|\mathcal{F}_i \langle T_i \rangle| \quad \mbox{for all $i \in [k]
\backslash \{j\}$,}\]
and that
\[|\mathcal{A}_j| \leq \left\{ \begin{array}{ll} {r \choose t}|
\mathcal{F}_j \langle T_j \rangle| & \mbox{if $\mathcal{A}_h \neq
\emptyset$ for some $h \in [k] \backslash \{j\}$;}\\
|\mathcal{F}_j| & \mbox{if $\mathcal{A}_i = \emptyset$ for all $i
\in [k] \backslash \{j\}$.}
\end{array} \right.\]
\medskip

\textit{Sub-case 1.1: $\mathcal{A}_i = \emptyset$ for all $i \in
[k] \backslash \{j\}$}. Then
\[\sum_{i=1}^k |\mathcal{A}_i| = |\mathcal{A}_j| \leq |\mathcal{F}_j| \leq \max\left\{\sum_{i = 1}^k |\mathcal{F}_i \langle T_i \rangle|, |\mathcal{F}_1|, \dots, |\mathcal{F}_k|\right\},\]
and equality holds throughout only if $|\mathcal{F}_j| = \max\left\{\sum_{i = 1}^k |\mathcal{F}_i \langle T_i \rangle|, |\mathcal{F}_1|, \dots, |\mathcal{F}_k|\right\}$ and $\mathcal{A}_j =
\mathcal{F}_j$, in which case either (ii) holds or (iii) holds.\medskip

\textit{Sub-case 1.2: $\mathcal{A}_h \neq \emptyset$ for some $h
\in [k] \backslash \{j\}$}. Then $|\mathcal{A}_j| \leq {r \choose
t}| \mathcal{F}_j \langle T_j \rangle|$. Let $x := \sum_{i \in [k]
\backslash \{j\}} |\mathcal{F}_i \langle T_i \rangle|$.

If $x \geq \frac{\mu - r}{\mu - r - (r-t){m \choose t+1}} {r
\choose t} |\mathcal{F}_j\langle T_j \rangle|$, then we have
\begin{align} &{r \choose t} |\mathcal{F}_j\langle T_j \rangle|
\leq \frac{\mu - r - (r-t){m \choose t+1}}{\mu - r}x \;
\Rightarrow \; |\mathcal{A}_j| \leq \left( 1 - \frac{r-t}{\mu - r}
{m \choose t+1} \right)x \nonumber \\
& \Rightarrow \; |\mathcal{A}_j| + \frac{r-t}{\mu - r} {m \choose
t+1}x \leq x \; \Rightarrow \; |\mathcal{A}_j| + \sum_{i \in [k]
\backslash \{j\}} \frac{r-t}{\mu - r} {m \choose t+1}
|\mathcal{F}_i \langle T_i \rangle| < \sum_{i = 1}^k
|\mathcal{F}_i \langle T_i \rangle| \nonumber \\
& \Rightarrow \; |\mathcal{A}_j| + \sum_{i \in [k] \backslash
\{j\}} |\mathcal{A}_i| < \sum_{i = 1}^k |\mathcal{F}_i \langle T_i
\rangle| \; \Rightarrow \; \sum_{i = 1}^k |\mathcal{A}_i| <
\sum_{i = 1}^k |\mathcal{F}_i \langle T_i \rangle|. \nonumber
\end{align}

Now suppose $x < \frac{\mu - r}{\mu - r - (r-t){m \choose t+1}} {r \choose t} |\mathcal{F}_j\langle T_j \rangle|$. By Lemma~\ref{maincor2}, $|\mathcal{F}_j| > \frac{{\mu - r \choose t}} {{r \choose t}} |\mathcal{F}_j \langle T_j \rangle|$. We will show that $\sum_{i=1}^k |\mathcal{A}_i| < |\mathcal{F}_j|$ follows from these last two inequalities together with those for $|\mathcal{A}_1|, \dots, |\mathcal{A}_k|$ above and the following. For any three integers $n, p$ and $l$ with $1 \leq l \leq p \leq n$,
\begin{equation}\left( \frac{n}{p} \right)^{l} \leq \prod_{q = 0}^{l-1}
\frac{n-q}{p-q} \leq \left( \frac{n-l+1}{p-l+1} \right)^l \nonumber %\label{binomin1}
\end{equation}
because, if $l \geq 2$, then $\frac{n}{p} \leq \frac{n-1}{p-1} \leq
\dots \leq \frac{n-l+1}{p-l+1}$; in particular, since ${n \choose
p} = \prod_{q = 0}^{p-1} \frac{n-q}{p-q}$,
\begin{equation} \left( \frac{n}{p} \right)^{p} \leq {n \choose p} \leq
(n-p+1)^p. \nonumber %\label{binomin2}
\end{equation}
We have
\begin{align} |\mathcal{F}_j| - \sum_{i=1}^k
|\mathcal{A}_i| &= |\mathcal{F}_j| - |\mathcal{A}_j| - \sum_{i
\in [k] \backslash \{j\}} |\mathcal{A}_i| \nonumber \\
&> \frac{{\mu - r \choose t}} {{r \choose t}} |\mathcal{F}_j \langle T_j \rangle| -
{r \choose t}| \mathcal{F}_j \langle T_j \rangle| - \sum_{i \in
[k] \backslash \{j\}} \frac{r-t}{\mu - r} {m \choose t+1}
|\mathcal{F}_i \langle T_i \rangle|
\nonumber \\
&= \left( \prod_{q = 0}^{t-1} \frac{(\mu-r)-q}{r-q} \right)
|\mathcal{F}_j \langle T_j \rangle| - {r \choose t}|
\mathcal{F}_j \langle T_j \rangle| - \frac{r-t}{\mu - r} {m
\choose t+1} x \nonumber \\
&> \left( \frac{\mu-r}{r} \right)^t |\mathcal{F}_j \langle T_j
\rangle| - {r \choose t}| \mathcal{F}_j \langle T_j \rangle| \nonumber \\
& \quad - \frac{r-t}{\mu - r} {m \choose t+1} \frac{\mu - r}{\mu - r -
(r-t){m \choose t+1}} {r \choose t} |\mathcal{F}_j\langle T_j
\rangle| \nonumber \\
&\geq |\mathcal{F}_j \langle T_j \rangle| \left( \left(
\frac{n_{\rm S} -r}{r} \right)^t - {r \choose t} - \frac{r - t}{n_{\rm S} - r - (r-t) {m \choose t+1}} {r \choose t} {m \choose t+1} \right)
\nonumber \\
&= |\mathcal{F}_j \langle T_j \rangle| \left( \left(
\frac{(r-t+1){m \choose t+1}}{r} \right)^t - {r \choose t} -
\frac{r - t}{{m \choose t+1}} {r \choose t} {m \choose t+1}
\right) \nonumber \\
&= |\mathcal{F}_j \langle T_j \rangle| \left( \frac{(r-t+1)^t{m
\choose t+1}^t}{r^t} - (r - t + 1){r \choose t} \right) \nonumber
\\
&\geq |\mathcal{F}_j \langle T_j \rangle| \left( \frac{(r-t+1)^t{2r
\choose t+1}^t}{r^t} - (r - t + 1)(r-t+1)^t \right) \nonumber \\
&\geq (r-t+1)^t|\mathcal{F}_j \langle T_j \rangle| \left(
\frac{1}{r^t} \left( \left( \frac{2r}{t+1} \right)^{t+1} \right)^t
- (r - t + 1) \right) \nonumber \\
&= (r-t+1)^t|\mathcal{F}_j \langle T_j \rangle| \left(
\frac{2^{t(t+1)}}{(t+1)^{t+1}} \left( \frac{r}{t+1}
\right)^{t^2-1} r - (r - t + 1) \right) \nonumber \\
&\geq (r-t+1)^t|\mathcal{F}_j \langle T_j \rangle| \left( \left(
\frac{2^t}{t+1} \right)^{t+1}r - (r - t + 1) \right) \quad
\mbox{(as $r \geq t+1$)}. \nonumber
\end{align}
By straightforward induction, $2^a \geq a+1$ for every integer $a \geq 1$. Thus, $2^t \geq t+1$, and hence $|\mathcal{F}_j| - \sum_{i=1}^k |\mathcal{A}_i| > (r-t+1)^t|\mathcal{F}_j \langle T_j \rangle| \left( r - (r - t + 1) \right) \geq 0$.

We have therefore shown that, for this sub-case, we have
$\sum_{i=1}^k |\mathcal{A}_i| < \sum_{i = 1}^k |\mathcal{F}_i
\langle T_i \rangle|$ or $\sum_{i=1}^k |\mathcal{A}_i| <
|\mathcal{F}_j|$, and hence $\sum_{i=1}^k |\mathcal{A}_i| <
\max\left\{\sum_{i = 1}^k |\mathcal{F}_i \langle T_i \rangle|,
|\mathcal{F}_1|, \dots, |\mathcal{F}_k|\right\}$.\medskip

\textit{Case 2: For all $i \in [k]$, $\mathcal{A}_i$ is empty or a trivial $t$-intersecting family.} Then, for all $i \in [k]$,
$\mathcal{A}_i$ is a subfamily of a $t$-star of $\mathcal{F}_i$,
and hence $|\mathcal{A}_i| \leq |\mathcal{F}_i \langle T_i
\rangle|$. Thus, $\sum_{i=1}^k |\mathcal{A}_i| \leq \sum_{i=1}^k
|\mathcal{F}_i \langle T_i \rangle|$.

Suppose $\sum_{i=1}^k |\mathcal{A}_i| = \sum_{i=1}^k
|\mathcal{F}_i \langle T_i \rangle|$. Then, for all $i \in [k]$,
$|\mathcal{A}_i| = |\mathcal{F}_i \langle T_i \rangle|$, and hence
$\mathcal{A}_i = \mathcal{F}_i \langle T_i' \rangle$ for some
$T_i' \in \mathcal{T}_i$. Let $j \in [2,k]$. If $S_j = S_1 = \{t\}$, then $\mathcal{A}_j = \{T_j'\}$, $\mathcal{A}_1 = \{T_1'\}$, and hence
$T_j' = T_1'$ (by the cross-$t$-intersection condition). Now suppose $S_j \neq \{t\}$ or $S_1 \neq \{t\}$. We have $r \geq t+1$, so $n_{\rm S} \geq 2r$.

Suppose $T_j' \neq T_1'$. Then, since $|T_j'| = |T_1'|$, $T_1'$
has an element $a_1$ that is not in $T_j'$. Let $p \in S_j$ and $q
\in S_1$. For each $l \in \{1,j\}$, let $M_l$ be a base of $\mathcal{H}$ such that $T_l \subseteq M_l$; we have $|M_l| \geq \mu(\mathcal{H}) \geq n_{\rm S} \geq 2r \geq p + q$. Thus, we can choose a $p$-subset $A'$ of $M_j \backslash \{a_1\}$ and a $q$-subset $B'$ of $(M_1 \backslash A') \cup T_1'$ such that $T_j' \subseteq A'$ and $T_1' \subseteq B'$.
Since $A' \subseteq M_j \in \mathcal{H}$ and $B' \subseteq M_1 \in
\mathcal{H}$, we have $A', B' \in \mathcal{H}$ (as $\mathcal{H}$
is hereditary). Thus, we now have $A' \in \mathcal{H}^{(p)} \langle
T_j' \rangle \subseteq \mathcal{F}_j \langle T_j' \rangle =
\mathcal{A}_j$ and $B' \in \mathcal{H}^{(q)} \langle T_1' \rangle
\subseteq \mathcal{F}_1 \langle T_1' \rangle = \mathcal{A}_1$.
However, $|A' \cap B'| = |A' \cap T_1'| = |A' \cap (T_1'
\backslash \{a_1\})| \leq t-1$, which contradicts the assumption
that $\mathcal{A}_j$ and $\mathcal{A}_1$ are
cross-$t$-intersecting.

Therefore, $T_j' = T_1'$. Since $j$ is an arbitrary element of
$[2,k]$, we get $T_1' = T_2' = \dots = T_k'$. Thus, we have
$T_1' \in \bigcap_{i=1}^k \mathcal{T}_i$ and $\mathcal{A}_i =
\mathcal{F}_i \langle T_1' \rangle$ for all $i \in [k]$.

We have shown that $\sum_{i = 1}^k |\mathcal{A}_i| \leq \max\left\{\sum_{i =
1}^k |\mathcal{F}_i \langle T_i \rangle|, |\mathcal{F}_1|, \dots, |\mathcal{F}_k|\right\}$, and that equality holds only if either (i) holds or (iii) holds.~\hfill{$\Box$} \\
\\
\textbf{Proof of Theorem~\ref{sumcor1}.} For each $i \in [k]$, let $T_i \in \mathcal{T}_i$. Let $\mu := \mu(\mathcal{H})$. In view of Theorem~\ref{sumresult}, it suffices to show that $\sum_{i = 1}^k |\mathcal{F}_i \langle T_i \rangle| < |\mathcal{F}_j|$. By Lemma~\ref{maincor2},
\begin{align} \sum_{i = 1}^k |\mathcal{F}_i \langle T_i \rangle|
&< \sum_{i = 1}^k \frac{{r \choose t}}{{\mu - r \choose t}}
|\mathcal{F}_i| \leq \sum_{i = 1}^k \frac{{r \choose t}}{{\mu - r
\choose t}} |\mathcal{F}_j| = k \left( \prod_{q = 0}^{t-1}
\frac{r-q}{(\mu-r)-q} \right) |\mathcal{F}_j| \nonumber \\
&\leq k \frac{r^t}{(\mu-r)^t} |\mathcal{F}_j| \leq \frac{kr^t}{(k^{1/t}r)^t} |\mathcal{F}_j| = |\mathcal{F}_j|,
\nonumber
\end{align}
as required.~\hfill{$\Box$}

\section{Proof of Theorem~\ref{prodresult}} \label{prodsection}

We now prove Theorem~\ref{prodresult}, which yielded all the other results in Section~\ref{prodsubsec}. %We now prove our main result for the product.
\\
\\
\textbf{Proof of Theorem~\ref{prodresult}.} As in Theorem~\ref{sumresult}, $\mathcal{T}_i \neq \emptyset$ for all $i \in [k]$. For each $i \in [k]$, let $T_i \in \mathcal{T}_i$. It is straightforward that $\prod_{i=1}^k |\mathcal{A}_i| = \prod_{i=1}^k |\mathcal{F}_i \langle T_i \rangle|$ if $\bigcap_{i = 1}^k \mathcal{T}_i \neq \emptyset$, $T \in \bigcap_{i = 1}^k \mathcal{T}_i$, and $\mathcal{A}_i = \mathcal{F}_i\langle T
\rangle$ for all $i \in [k]$. We now show that $\prod_{i = 1}^k
|\mathcal{A}_i| \leq \prod_{i = 1}^k |\mathcal{F}_i \langle T_i
\rangle|$, and that equality holds only if $\bigcap_{i = 1}^k \mathcal{T}_i \neq \emptyset$ and, for some $T \in \bigcap_{i = 1}^k
\mathcal{T}_i$, $\mathcal{A}_i = \mathcal{F}_i\langle T \rangle$
for all $i \in [k]$.

If one of $\mathcal{A}_1, \dots, \mathcal{A}_k$ is empty, then $\prod_{i = 1}^k |\mathcal{A}_i| = 0 < \prod_{i = 1}^k |\mathcal{F}_i \langle T_i \rangle|$. Now suppose $\mathcal{A}_i \neq \emptyset$ for all $i \in [k]$.\medskip

\textit{Case 1: For some $j \in [k]$, $\mathcal{A}_j$ is not a
trivial $t$-intersecting family.} By the argument for the corresponding case in the proof of Theorem~\ref{sumresult}, we obtain
\[|\mathcal{A}_j| \leq {r \choose t} |\mathcal{F}_j \langle T_j \rangle|\]
(as $\mathcal{A}_i \neq \emptyset$ for all $i \in [k]$) and
\[|\mathcal{A}_i| < \frac{r-t}{\mu(r,t) - r} {m(r,t) \choose t+1}
|\mathcal{F}_i \langle T_i \rangle| \quad \mbox{for all $i \in [k]
\backslash \{j\}$.}\]
We therefore have
\begin{align}\prod_{i = 1}^k |\mathcal{A}_i| &= |\mathcal{A}_j|
\prod_{i \in [k] \backslash \{j\}} |\mathcal{A}_i| < {r \choose
t} |\mathcal{F}_j \langle T_j \rangle| \prod_{i \in [k]
\backslash \{j\}} \left( \frac{r-t}{n_{\rm P}(r,t) - r} {m(r,t)
\choose t+1} |\mathcal{F}_i \langle T_i \rangle| \right) \nonumber \\
&= {r \choose t} \left( \frac{r-t}{n_{\rm P}(r,t) - r} {m(r,t) \choose
t+1}\right)^{k-1} \prod_{i=1}^k |\mathcal{F}_i \langle T_i
\rangle| = {r \choose t} \left( \frac{1}{{r \choose t}}
\right)^{k-1} \prod_{i = 1}^k |\mathcal{F}_i \langle T_i \rangle|,
\nonumber
\end{align}
and hence $\prod_{i = 1}^k |\mathcal{A}_i| < \prod_{i = 1}^k
|\mathcal{F}_i \langle T_i \rangle|$.\medskip

\textit{Case 2: For all $i \in [k]$, $\mathcal{A}_i$ is a trivial
$t$-intersecting family.} An argument similar to that
for the corresponding case in the proof of
Theorem~\ref{sumresult} gives us that $\prod_{i = 1}^k |\mathcal{A}_i|
\leq \prod_{i = 1}^k |\mathcal{F}_i \langle T_i \rangle|$, and that
equality holds only if $\bigcap_{i = 1}^k \mathcal{T}_i \neq \emptyset$
and, for some $T \in \bigcap_{i = 1}^k \mathcal{T}_i$,
$\mathcal{A}_i = \mathcal{F}_i\langle T \rangle$ for all $i \in
[k]$.~\hfill{$\Box$}

\end{document}